\documentclass[10pt]{amsart}
\usepackage{amsmath,amscd,latexsym,verbatim,amssymb}
\usepackage[francais]{babel}

  \newcommand{\resp}{{\it resp.} }

  \newcommand{\Z}{\mathbf{Z}}

\theoremstyle{definition}

\newcounter{spec}
\begin{document}
  \centerline{\bf Id\'ees galoisiennes}
 
  \medskip  
 \centerline{Yves Andr\'e}

      \bigskip 
  
  \bigskip
  
       \bigskip
       
  \subsection{}    On c\'el\`ebre cette ann\'ee le bicentenaire d'\'Evariste Galois, n\'e le 25 octobre 1811.

 Galois est \`a l'origine de la notion de groupe, cette notion fondamentale qui allait, sous divers avatars, envahir toutes les math\'ematiques, ainsi qu'une bonne part de la physique et m\^eme de la chimie. Dans ses R\'ecoltes et Semailles, A. Grothendieck n'h\'esite pas \`a parler de l'invention du z\'ero et de l'id\'ee de groupe et comme des deux plus grandes innovations math\'ematiques de tous les temps.

Galois est aussi le premier \`a avoir formul\'e le principe de correspondance entre sym\'etries et invariants, qui s'av\'erera f\'econd, comme il semble l'avoir lui-m\^eme pressenti, bien au-del\`a du contexte originel de la th\'eorie des \'equations alg\'ebriques. Ê

ÊCes id\'ees sem\'ees, qui en germant ont r\'evolutionn\'e les math\'ematiques, ont \'et\'e consign\'ees sur quelques dizaines de feuillets par un jeune homme qui mourut dans sa vingti\`eme ann\'ee. 
 La veille du duel fatal, il \'ecrivit une splendide lettre-testament, qui d\'ebute ainsi:
 
 \begin{quotation}{\small   JÕ'ai fait en analyse plusieurs choses nouvelles.
Les unes concernent la th\'eorie des \'equations, les autres les fonctions int\'egrales.
Dans la th\'eorie des \'equations, j'ai recherch\'e dans quels cas les \'equations \'etaient r\'esolubles par des radicaux : ce qui m'a donn\'e occasion d'approfondir cette th\'eorie, et de d\'ecrire toutes les transformations possibles sur une \'equation lors m\^eme qu'elle n'est pas soluble par radicaux.}\end{quotation}

C'est dans cette lettre qu'il appelle {\og}{th\'eorie de l'ambigu¬{\i}t\'e}{\fg} le corpus d'id\'ees \'elabor\'e par lui pour \'elucider la th\'eorie des \'equations alg\'ebriques.

  \subsection{}   Nous proposons un libre parcours autour de l'h\'eritage de ces id\'ees. 
  
  Les domaines les plus classiques d'application de la {\og}{th\'eorie de l'ambigu¬{\i}t\'e}{\fg} concernent les {\it nombres alg\'ebriques} et les {\it fonctions alg\'ebriques}.
Nous les survolerons, avant d'aborder ceux des {\it fonctions transcendantes} et des {\it nombres transcendants} qui semblent avoir aussi \'et\'e envisag\'es par Galois, comme en t\'emoigne la fin de sa lettre:

   \begin{quotation}{\small Tu sais, mon cher Auguste, que ces sujets ne sont pas les seuls que j'ai explor\'es. Mes principales m\'editations depuis quelque temps \'etaient dirig\'ees sur l'application \`a l'analyse transcendante de la th\'eorie de l'ambigu\^{\i}t\'e. Il s'agissait de voir a priori dans une relation entre quantit\'es ou fonctions transcendantes quels \'echanges on pouvait faire, quelles quantit\'es on pouvait substituer aux quantit\'es donn\'ees sans que la relation p\^ut cesser d'avoir lieu. Cela fait reconna\^{\i}tre tout de suite l'impossibilit\'e de beaucoup d'expressions que l'on pouvait chercher. Mais je n'ai pas le temps et mes id\'ees ne sont pas encore bien d\'evelopp\'ees sur ce terrain qui est immense.}\end{quotation}

      \section{Nombres alg\'ebriques et groupes de Galois}

  \subsection{}  Dans la premi\`ere encyclop\'edie math\'ematique de la renaissance\footnote{Summa de arithmetica, geometria, de proportioni et de proportionalita (Venise, 1494).}, l'auteur-compilateur L. Pacioli compare l'impossibilit\'e de r\'esoudre les \'equations de degr\'e 3 \`a l'impossibilit\'e de la quadrature du cercle.
 Pourtant, une vingtaine d'ann\'ees plus tard, son coll\`egue S. Del Ferro trouvait la formule par radicaux. 
 On conna\^it la saga de secrets, querelles et trahisons qui s'ensuivit. A la fin XVIe, Bombelli \'etait en mesure de donner un expos\'e complet de la r\'esolution par radicaux des \'equations de degr\'e au plus 4, y compris une discussion des ambiguit\'es et des nombres imaginaires comme $\sqrt -1$ qui interviennent dans les formules.
 
La r\'esolution par radicaux des \'equations d'ordre au moins 5 a r\'esist\'e... jusqu'\`a ce qu'Abel en d\'emontre l'impossibilit\'e en g\'en\'eral, et que Galois fasse toute la lumi\`ere sur cette question\footnote{ Galois avait d\'emontr\'e, \`a dix-sept ans, cette impossibilit\'e; avant de d\'ecouvrir qu'Abel venait de publier le r\'esultat (Beweis der Unmšglichkeit der algebraischen Auflšsbarkeit der allgemeinen Gleichungen, J. de Crelle, 1826, 65-84), et \'ecrivit une courte note exposant son approche. Abel mourut de tuberculose, deux ans avant Galois, \`a l'\^age de 26 ans.}.
 
 \smallskip Les progr\`es sont venus d'une r\'eflexion de plus en plus pouss\'ee autour des sym\'etries entre racines d'un polyn\^ome (Vi\`ete, van der Monde, Leibniz, Lagrange, Cauchy, Abel) qui culmine dans la th\'eorie de Galois. 

  \subsection{}   Avec Galois, les ambigu¬{\i}t\'es ne constituent plus une nuisance, elles constituent un groupe!  
       
  Soit $\alpha$ un {\it nombre alg\'ebrique}, racine d'un polyn\^ome $  x^n + a_{n-1}x^{n-1}+ \cdots + a_1 x + a_0  $  \`a coefficients rationnels, de degr\'e $n$ minimal. Factorisons-le en $(x-\alpha_1)\cdots (x-\alpha_n) $, o\`u les $\alpha_i$ (les autres racines) sont les {\it conjugu\'es} de $\alpha=\alpha_1$. 
   
      Les $\alpha_i$ engendrent un {\og}{corps de nombres}{\fg}, et le {\it groupe de Galois} associ\'e \`a $\alpha$  est le groupe d'automorphismes de ce corps (c'est-\`a-dire le groupe des permutations des $\alpha_i$ qui respectent l'addition et la multiplication)\footnote{Galois avait aussi imagin\'e de prendre pour $a_i$ des entiers modulo un nombre premier, et avait construit ainsi les {\og}{corps finis}{\fg} (imaginaires de Galois).}.

   \smallskip   Le point de vue de Galois lui-m\^eme \'etait sensiblement diff\'erent (un si\`ecle de r\'e\'elaborations s\'eparent ces points de vue\footnote{on renvoie \`a la th\`ese de C. Ehrhardt (Paris, 2007) pour une analyse d\'etaill\'ee de ces r\'e\'elaborations. La notion m\^eme de groupe a elle aussi subi de nombreuses vicissitudes, la condition d'associativit\'e n'ayant \'et\'e clairement d\'egag\'ee qu'au d\'ebut du XX\`eme si\`ecle.}).
    Il est bien rendu par cette analyse de G. Ch\^atelet:

                 \begin{quotation}{\small La th\'eorie de Galois constitute probablement un des plus beaux exemples du principe de dissym\'etrie cr\'eatrice en math\'ematiques. [...] Ces racines passeront dans l'existence concr\`ete dans le domaine $D(\alpha_1, \cdots, \alpha_n)$ lorsque je ne serai pas seulement capable d'\'ecrire un signe {\og}{formel}{\fg} mais capable aussi d'exhiber un proc\'ed\'e de discernement entre ces racines. Ce proc\`es de discernement est bien une cassure de la sym\'etrie du groupe des racines. Discerner plus, c'est se montrer capable d'exhiber certaines {\og}{grandeurs tests}{\fg}, non rationnelles, et invariantes par un groupe de sym\'etrie plus restreint. Les racines vont s'individuer au fur et \`a mesure de la pr\'ecision des {\og}{expressions alg\'ebriques test}{\fg} construites par certaines adjonctions aux rationnels. [...]
        
        {\it Rendre raison} de l'existence conscr\`ete de ces racines est \'equivalent \`a la {\it pr\'esentation du groupe de sym\'etrie} de l'\'equation comme une cha\^{\i}ne descendante construite de telle mani\`ere que l'\'ecart entre deux de ses maillons soit le manque de discernement  relatif entre les racines lorsque je calcule dans le domaine associ\'e.\footnote{G. Ch\^atelet, La physique math\'ematique comme projet, L'enchantement du virtuel, p. 115-117, Presses de l'ENS}} \end{quotation}

       La correspondance de Galois met en regard extensions de corps obtenues par adjonction de racines, et certains groupes de sym\'etrie de ces racines. Elle remplace, si on veut, l'imbrication de deux op\'erations commutatives (l'addition et la multiplication) dans des extensions, par une seule op\'eration non commutative (la composition des permutations). 
     
       \subsection{} Lorsque l'on adjoint de plus en plus de nombres alg\'ebriques, la correspondance de Galois donne lieu \`a un syst\`eme {\og}{projectif}{\fg} de groupes finis s'envoyant les uns sur les autres, et dont la {\og}{limite}{\fg} ${\rm{Gal}}({\bar{\mathbb Q}}/{\mathbb Q})$, 
   appel\'ee {\it groupe de Galois absolu de ${\mathbb Q}$}, n'est autre que le groupe infini\footnote{naturellement muni d'une structure de groupe topologique compact.}
 des automorphismes du corps ${\bar{\mathbb Q}}$ des nombres (complexes) alg\'ebriques.  
       
           \smallskip  C'est cette d\'emarche qu'A. Lautman appelle la {\og}{mont\'ee vers l'absolu}{\fg}\footnote{Essai sur les notions de structure et d'existence en math\'ematiques. R\'e\'edition Vrin 2006.}: un seul objet math\'ematique, ${\rm{Gal}}({\bar{\mathbb Q}}/{\mathbb Q})$, code les propri\'et\'es de toutes les \'equations alg\'ebriques \`a coefficients rationnels \`a la fois. 
Ce groupe de Galois absolu, objet central de la th\'eorie des nombres, reste largement un myst\`ere apr\`es un si\`ecle et demi d'efforts intenses pour en comprendre la structure.   

Ses quotients commutatifs sont d\'ecrits par la cyclotomie, qui remonte \`a Gauss. Pour aller au-del\`a, on peut faire agir lin\'eairement ce groupe, t\^acher d'\'etudier ses repr\'esentations lin\'eaires\footnote{comme l'a montr\'e bien plus tard Grothendieck, la cohomologie \'etale des vari\'et\'es alg\'ebriques fournit de telles repr\'esentations galoisiennes, tandis que le programme de Langlands relie certaines repr\'esentations galoisiennes \`a certains objets de l'analyse spectrale: les repr\'esentations automorphes.}. 
  La th\'eorie des repr\'esentations lin\'eaires est d'ailleurs n\'ee, entre les mains de Frobenius, \`a propos d'un probl\`eme que lui avait soumis Dedekind en 1896 concernant des groupes de Galois de corps de nombres. Bien des conqu\^etes r\'ecentes de la th\'eorie des nombres (et notamment la preuve par A. Wiles du {\og}{th\'eor\`eme{\fg}{ de Fermat) reposent sur l'\'etude des d\'eformations de ces objets {\og}{souples}{\fg} que sont les repr\'esentations galoisiennes. 

 Une autre fa\c con de {\og}{comprendre}{\fg} ${\rm{Gal}}({\bar{\mathbb Q}}/{\mathbb Q})$ est d'en d\'ecrire les quotients finis: le fameux probl\`eme de Galois inverse, toujours ouvert, consiste \`a savoir si tout groupe fini est groupe de Galois d'un corps de nombres. 

Une troisi\`eme fa\c con, apparemment beaucoup plus modeste, consisterait \`a distinguer, nommer ou construire des \'el\'ements particuliers. Or il faut avouer qu'apr\`es l'identit\'e et la conjugaison complexe, on ne sait pas d\'ecrire un troisi\`eme \'el\'ement, m\^eme si la th\'eorie des {\og}{dessins d'enfants}{\fg} de Grothendieck vise une compr\'ehension graphique des groupes de Galois et irait dans ce sens\footnote{on renvoie \`a l'expos\'e de A. Zvonkine aux Journ\'ees math\'ematiques X-UPS 2004 pour une introduction \`a cette th\'eorie.}.

    \section{Fonctions alg\'ebriques et groupes de Galois} 
 
 \subsection{}  Galois \'etait bien conscient que sa {\og}{th\'eorie de l'amgu\"{\i}t\'e s'appliquait aussi bien aux fonctions alg\'ebriques qu'aux nombres alg\'ebriques. Cela ressort de son calcul du groupe de Galois des points de division d'une courbe elliptique, et de ses myst\'erieux travaux sur les int\'egrales ab\'eliennes\footnote{on renvoie \`a l'expos\'e de P. Popescu dans ce volume pour une introduction aux int\'egrales ab\'eliennes et une discussion des ambigu\"{\i}t\'es qu'elles charrient. On ignore comment Galois avait pu obtenir les r\'esultats  sur les int\'egrales ab\'eliennes expos\'es dans sa lettre, qui pr\'efigurent ceux de B. Riemann un quart de si\`ecle plus tard.}, qu'il introduit ainsi dans sa lettre:
 
 \begin{quotation}{\small
   On traite \`a la fois toutes les int\'egrales dont la diff\'erentielle est une fonction de la variable et d'une m\^eme fonction irrationnelle de la variable, que cette irrationnelle soit ou ne soit pas un radical, qu'elle s'exprime ou ne s'exprime pas par des radicaux.
 }\end{quotation}
 
   Mais c'est F. Klein\footnote{l'histoire commence en 1877 lorsque Klein remarque que le groupe des isom\'etries laissant invariant l'icosa\`edre est isomorphe au groupe de Galois d'une \'equation quintique \`a coefficients dans le corps de fonctions rationnelles d'une variable auxiliaire.} qui en croisant la th\'eorie de l'ambiguit\'e et celle des surfaces \`a plusieurs feuillets de Riemann, a inaugur\'e une tradition qui, via la th\'eorie de l'uniformisation et du groupe fondamental de Poincar\'e, aboutira \`a la fusion qu'A. Grothendieck effectuera vers 1960 entre th\'eorie de Galois et th\'eorie des rev\^etements, pierre de touche de la fusion plus g\'en\'erale qu'il a mise en oeuvre entre g\'eom\'etrie alg\'ebrique et th\'eorie des nombres.

 \subsection{}  En langage moderne, voici de quoi il s'agit. Une fonction alg\'ebrique est racine d'un polyn\^ome \`a coefficients fonctions rationnelles d'une variable auxiliaire $t$.  Les solutions complexes de l'\'equation  forment  une surface de Riemann, qui est un rev\^etement du plan complexe (muni de la coordonn\'ee complexe $t$). 
  
\smallskip Plus g\'en\'eralement, on a la notion de rev\^etement fini normal 
$Y\to X$ de surfaces de Riemann (qu'on suppose non ramifi\'e pour simplifier) . Dans ce contexte g\'eom\'etrique, ce qui remplace les racines d'une \'equation, ce sont les points $\{y_1,\ldots , y_n\}$   de $Y$ qui s'envoient sur un point $x$ arbitraire fix\'e de $X$. Le groupe de Galois ${\rm{Gal}}(Y/X)$ est un sous-groupe du groupe des permutations de $\{y_1,\ldots , y_n\}$\footnote{dans le cas de la quintique de Klein, le groupe de Galois est le groupe de l'icosa\`edre, auquel Klein a consacr\'e un ouvrage classique.}. Il peut se calculer comme suit: tra\c cons sur $X$ un lacet $\gamma$ point\'e en $x$, et choisissons un point $y_i$ de $Y$ au-dessus de $x$. Un tel $\gamma$ se {\og}{rel\`eve}{\fg} alors en un chemin sur $Y$ partant de $y_i$, qui en g\'en\'eral aboutira \`a un autre point $y_j$, d'o\`u une permutation $y_i\mapsto y_j$ de l'ensemble des points au-dessus de $x$, qui ne d\'epend en fait du lacet $\gamma$ qu'\`a d\'eformation (homotopie) pr\`es. C'est ainsi que s'obtiennent les \'el\'ements de ${\rm{Gal}}(Y/X)$.

 \subsection{}  Dans ce contexte, on peut encore effectuer une {\og}{mont\'ee vers l'absolu}{\fg}. Ce qui correspond au corps ${\bar{\mathbb Q}}$ est maintenant le {\it rev\^etement universel} de $X$ (point\'e en $x$), et son groupe d'automorphismes n'est autre que le {\it groupe fondamental} de Poincar\'e $\pi_1(X, x)$: c'est groupe des lacets trac\'es sur $X$ \`a homotopie pr\`es, partant et aboutissant \`a un point $x$ fix\'e. 

\smallskip  On peut alg\'ebriser la construction en rempla\c cant $\pi_1(X, x )$ par la limite projective $\hat\pi_1(X, x )$ des groupes ${\rm{Gal}}(Y/X)$, lorsque le rev\^etement $Y/X$ grossit. Dand le point de vue cat\'egorique de Grothendieck, $\hat\pi_1(X, x )$ s'interpr\`ete comme groupe des automorphismes du foncteur fibre en $x$ $$ Y \mapsto Y_x =  \{y_1,\ldots , y_n\}$$  sur la cat\'egorie des rev\^etements de $(X, x)$ (\`a valeurs dans la cat\'egorie des ensembles), ce qui permet d'unifier les th\'eories de Galois arithm\'etique et g\'eom\'etrique dans un m\^eme moule.

 \subsection{}  Une approche indirecte fascinante de ${\rm{Gal}}({\bar{\mathbb Q}}/{\mathbb Q})$ consiste \`a le relier \`a la g\'eom\'etrie des surfaces de Riemann.  

Pour fixer les id\'ees, prenons pour $X$ le plan complexe priv\'e des points $0$ et $1$. Son groupe fondamental $\pi_1(X )$ est le groupe libre \`a deux g\'en\'erateurs, les lacets $\gamma_0$ et $\gamma_1$ autour de $0$ et de $1$. Son compl\'et\'e $\hat\pi_1(X )$ s'obtient en permettant des {\og}{mots infinis}{\fg} en $\gamma_0$ et $\gamma_1$ (et leurs inverses). 

Suivant Grothendieck et Belyi, ${\rm{Gal}}({\bar{\mathbb Q}}/{\mathbb Q})$ op\`ere sur les rev\^etements finis de $X$, donc sur $\hat\pi_1(X )$, et cette op\'eration est {\it fid\`ele: le groupe de Galois absolu arithm\'etique ${\rm{Gal}}({\bar{\mathbb Q}}/{\mathbb Q})$ se plonge dans le groupe des automorphismes du groupe de Galois absolu g\'eom\'etrique $\hat\pi_1(X )$}.  
 
 De l\`a, Grothendieck a alors propos\'e de d\'ecrire ${\rm{Gal}}({\bar{\mathbb Q}}/{\mathbb Q})$ au moyen de notions graphiques {\og}{si simples qu'un enfant peut les conna\^{\i}tre en jouant}{\fg}. Consid\'erons un rev\^etement $Y\to X$, en supposant pour simplifier que $Y$ est le plan complexe priv\'e de quelques points. L'image inverse dans $Y$ du segment $]0, 1[$ de $X$ est un objet combinatoire tr\`es simple que Grothendieck appelle {\og}{dessin d'enfant}{\fg}. Le d\'efi est de comprendre en termes combinatoires l'op\'eration fid\`ele de ${\rm{Gal}}({\bar{\mathbb Q}}/{\mathbb Q})$ sur ces dessins\footnote{sur tout cela, on pourra consulter R. Douady, A. Douady, Alg\`ebgre et th\'eories galoisiennes, Cassini 2005.}. 
. 
 
 Pour ce faire, il faut disposer au pr\'ealable d'un codage combinatoire des \'el\'ements de ${\rm{Gal}}({\bar{\mathbb Q}}/{\mathbb Q})$ lui-m\^eme. C'est ce qui a \'et\'e obtenu par V. Drinfeld autour de 1990 (en d\'ecouvrant un lien insoup\c conn\'e entre ${\rm{Gal}}({\bar{\mathbb Q}}/{\mathbb Q})$ et groupes quantiques): il plonge ${\rm{Gal}}({\bar{\mathbb Q}}/{\mathbb Q})$ dans un groupe de nature combinatoire ${\rm{GT}}$ - le groupe de Grothendieck-Teichm\"uller, d\'efini par g\'en\'erateurs (topologiques) et trois relations tr\`es simples -, qui s'av\`ere agir fid\`element sur les dessins d'enfants.  On ignore \`a l'heure actuelle si ${\rm{GT}}$ est r\'eellement {\og}{plus gros}{\fg} que ${\rm{Gal}}({\bar{\mathbb Q}}/{\mathbb Q})$ (on le soup\c conne).   

     \section{Fonctions transcendantes et groupes de Galois} 
     
  \subsection{}   J. Liouville, qui a exhum\'e,  \'etudi\'e et fait conna\^{\i}tre les papiers de Galois, est aussi le premier \`a avoir poursuivi les intuitions galoisiennes vers les fonctions transcendantes: au lieu de se demander quand une \'equation alg\'ebrique est r\'esoluble par radicaux, il se demande quand une \'equation diff\'erentielle lin\'eaire est r\'esoluble par quadratures (int\'egrales et exponentielles d'int\'egrales). Mais c'est E. Picard qui, en 1883, sous doute inspir\'e par des id\'ees de S. Lie, a introduit dans ce contexte le {\it groupe de Galois diff\'erentiel}: c'est le groupe form\'e des automorphismes commutant \`a la d\'erivation, parmi tous les automorphismes de l'extension du corps de fonctions de base obtenue en adjoignant les solutions de l'\'equation diff\'erentielle et leurs d\'eriv\'ees.
    
 Du fait que les solutions d'une \'equation diff\'erentielle lin\'eaire forment non plus un ensemble fini, mais un espace vectoriel de dimension finie, le groupe de Galois diff\'erentiel n'est plus un groupe fini en g\'en\'eral, mais un groupe alg\'ebrique (groupe de matrices). 
 C'est d'ailleurs, via les travaux d'E. Kolchin, l'une des sources historiques de la th\'eorie des groupes alg\'ebriques. 
 
 On dispose d'une correspondance de Galois diff\'erentielle mettant en regard extensions de corps diff\'erentiels obtenues par adjonction de solutions, et certains groupes de transformations lin\'eaires des espaces de solutions.
 
  \subsection{}  La th\'eorie a m\^uri lentement. La classification des ambigu\"{\i}t\'es galoisiennes dans le cadre des \'equations diff\'erentielles lin\'eaires analytiques au voisinage d'une singularit\'e, due \`a J. P. Ramis, date seulement de la fin du XX$^e$ si\`ecle. Le r\'esultat est que {\it dans ce cadre analytique local, il y a trois types, et trois seulement, d'ambigu\"{\i}t\'es galoisiennes} (qui sont des \'el\'ements du groupe de Galois diff\'erentiel):

\smallskip 1) la {\it monodromie}: c'est l'ambigu\"{\i}t\'e qui r\'esulte de ce que l'on ne retombe pas la valeur initiale lorsque l'on fait subir \`a une solution un tour autour de la singularit\'e\footnote{cette notion est due \`a Riemann, qui l'avait mise en \'evidence dans le contexte des \'equations diff\'erentielles hyperg\'eom\'etriques. Mais Galois a probablement pu avoir l'intuition de cette ambigu\"{\i}t\'e, notamment \`a propos des p\'eriodes d'int\'egrales ab\'eliennes d\'ependant d'un param\`etre, comme il ressort du passage de sa lettre o\`u il parle des {\og}{p\'eriodes relatives \`a une m\^eme {\it r\'evolution}}{\fg}. 
  }.  Consid\'erons par exemple l'\'equation diff\'erentielle 
$$ y' = \frac{1}{2x} y$$ au voisinage de la singularit\'e $0$; une solution est $y= \sqrt x$, et un tour autour de l'origine la transforme en $- y$. Le groupe de Galois diff\'erentiel de cette \'equation est le groupe \`a deux \'el\'ements engendr\'e par la monodromie.

\smallskip 2) le {\it recalibrage des exponentielles}: consid\'erons l'\'equation diff\'erentielle 
$$ xy' + y= 0$$ au voisinage de la singularit\'e $0$; une solution est $y= e^{1/x}$, et toute autre solution non nulle s'obtient en multipliant $y$ par une constante non nulle. Le groupe de Galois diff\'erentiel de cette \'equation est le groupe multiplicatif ${ \mathbb C}^\times$ engendr\'e par ces recalibrages.

\smallskip 3) les {\it ambigu\"{\i}t\'es de Stokes}: consid\'erons l'\'equation diff\'erentielle inhomog\`ene
$$ xy' + y= x$$ au voisinage de la singularit\'e $0$. L. Euler l'avait d\'ej\`a rencontr\'ee dans son fameux m\'emoire sur les s\'eries divergentes\footnote{De seriebus divergentibus, publi\'e en 1760 \`a l'Acad\'emie de St. Petersbourg. Voir aussi http://www.maa.org/news/howeulerdidit.html, juin 2006.

Euler n'h\'esitait pas \`a braver la divergence en \'ecrivant des formules comme:
    $$1+2+4+8+ 16+ \cdots = -1,\;\; \;\; 1+2+3+4+5+ \cdots = -1/12,$$ 
 ou  en attribuant une valeur pr\'ecise \`a la somme 
   $\; 1- 1! + 2! -3! +4!  -  5! + \cdots  \;$ (qu'il \'evalue num\'eriquement, de six mani\`eres diff\'erentes). 
Ces formules {\og}{scandaleuses}{\fg}, s\'ev\`erement critiqu\'ees aux temps de la qu\^ete de la rigueur en Analyse, furent pleinement \'eclaircies et justifi\'ees ult\'erieurement (\`a l'abus de notation pr\`es qu'elles commettent). Par exemple, la premi\`ere formule n'est autre que l'\'evaluation en $x=1$ de la s\'erie de puissances 
 $ \,  1+ 2x + 4x^2 + \cdots = 1/(1-2x)  $; stricto sensu, c'est la valeur en $1$ du prolongement analytique de la fonction $1/(1-2x)  $ de la variable complexe $x$ d\'efinie par cette s\'erie. 
  La seconde formule est nettement plus profonde et attendit 120 ans sa justification: elle exprime la valeur en $s=-1$ du prolongement analytique de la fonction  $ \zeta(s)   = \sum_1^\infty\, n^{-s} = \prod_p\, \frac{ 1}{1- p^{-s}}  $  de Riemann. Dans son article visionnaire  de 1859 (\"Uber die Anzahl der Primzahlen unter einer gegebenen Gr\"osse,  Monatsberichte der Berliner Akademie), Riemann prouve, pour tout nombre complexe $s$, la sym\'etrie sugg\'er\'ee par Euler entre $\zeta(s)$ et $\zeta(1-s)$, et il explique le lien entre la distribution des nombres premiers et la position des z\'eros de $\zeta$. }: c'est l'\'equation satisfaite par la s\'erie formelle $\hat y= \sum (-1)^n n! x^{n+1}$, qui diverge en tout point $x\neq 0$. Euler utilisait d'ailleurs cette \'equation pour  {\og}{sommer}{\fg} cette s\'erie divergente, en identifiant la {\og}{somme}{\fg} \`a la {\og}{vraie}{\fg} solution $y= \int_0^\infty \,\frac{e^{-t/x}}{1+t}dt$, dont $\hat y$ est le d\'eveloppement asymptotique (Euler \'ecrit \og{evolutio}\fg) dans le plan priv\'e de la demi-droite r\'eelle n\'egative.

Mais dans un autre secteur, le d\'eveloppement asymptotique de $y$ peut changer. Ces ambiguit\'es li\'ees au choix des secteurs donnent lieu \`a des ambigu\"{\i}t\'es galoisiennes, les ambigu\"{\i}t\'es de 
Stokes\footnote{G. Stokes les a mises en \'evidence, au milieu du XIX\`eme si\`ecle, sur l'\'equation diff\'erentielle d'Airy $y''= xy$, apr\`es avoir remarqu\'e qu'il \'etait bien plus efficace de calculer les z\'eros de la solution d'Airy en se servant du d\'eveloppement asymptotique divergent, \`a l'infini, plut\^ot qu'avec le d\'eveloppement de Taylor convergent \`a l'origine comme faisait G. Airy.}.

  \smallskip De mani\`ere g\'en\'erale, le groupe de Galois diff\'erentiel est engendr\'e (au sens des groupes alg\'ebriques) par ces trois types de matrices\footnote{voir par exemple M. van der Put, M. Singer, Galois theory of linear differential equations, Springer Grundlehren der Math. Wiss. 328, 2003.}.

    \section{Nombres transcendants et groupes de Galois}

    \subsection{} 
\medskip Peut-on attacher \`a un nombre transcendant  donn\'e des conjugu\'es, et un groupe de Galois qui les permute, comme dans le cas des nombres alg\'ebriques? 

\smallskip Commen\c cons par consid\'erer le cas de $\pi$, dont F. von Lindemann a montr\'e en 1882 qu'il est transcendant, c'est-\`a-dire ne satisfait aucune \'equation polyn\^omiale \`a coefficients rationnels (d\'emontrant ipso facto l'impossibilit\'e de la quadrature du cercle). Mais, comme Euler l'avait observ\'e, $\pi$ satisfait une telle \'equation mais de degr\'e infini (o\`u une s\'erie de puissances remplace un polyn\^ome)\footnote{formellement, l'\'equation d'Euler $\sum n^{-2}= \pi^2/6$ r\'esulte de l'extension \`a cette s\'erie de la formule de Newton pour la somme des carr\'es des racines en terme des coefficients d'une \'equation polyn\^omiale. Galois avait d'ailleurs fait la remarque que ce type de formules ne d\'ependait pas du nombre de racines, ayant sans doute en vue la possibilit\'e de faire tendre ce nombre vers l'infini...}
$$  \prod_{n\in \Z\setminus 0}\, (1-\frac{x}{n\pi}) = \frac{\sin x}{x} = 1 - \frac{x^2}{6} + \frac{x^4}{120} + \cdots \in {\mathbb Q}[[x]].$$
Ce qui sugg\`ere de consid\'erer les multiples de $\pi$ comme ses conjugu\'es. En fait, si l'on veut qu'un groupe permute transitivement les conjugu\'es, on est amen\'e \`a prendre pour conjugu\'es de $\pi$ tous ses multiples rationnels non nuls - le groupe de Galois serait alors le groupe multiplicatif $\mathbb Q^\times$. 

 \smallskip Peut-on g\'en\'eraliser cette approche?
 Un vieux r\'esultat peu connu de A. Hurwitz assure que tout nombre complexe est racine d'une s\'erie de puissances \`a coefficients rationnels qui converge partout. Peut-on alors consid\'erer ses autres racines comme des conjugu\'es? 
 
 H\'elas, cette approche est un cul-de-sac\footnote{sauf \`a imposer de fortes contraintes arithm\'etiques sur les d\'enominateurs des coefficients. Cette possibilit\'e est d'ailleurs loin d'avoir \'et\'e explor\'ee syst\'ematiquement.}, car il existe une infinit\'e ind\'enombrable de telles s\'eries, et aucun moyen d'en choisir une canonique en g\'en\'eral.
 
   \subsection{}  Nous allons voir qu'on peut tout de m\^eme s'attendre \`a pouvoir d\'efinir des conjugu\'es et un groupe de Galois pour une vaste classe de nombres (en g\'en\'eral transcendants), incluant la plupart des constantes math\'ematiques classiques. Ces nombres sont des int\'egrales,  les {\it p\'eriodes}\footnote{le nom vient de ce que, dans le cas particulier des courbes alg\'ebriques (d\'efinies sur un corps de nombres), ce sont les p\'eriodes d'int\'egrales ab\'eliennes, que consid\'erait Galois dans le passage \'evoqu\'e plus haut o\`u il parle de \og{r\'evolution}\fg, et qui apparaissent comme p\'eriodes (au sens usuel) des fonctions ab\'eliennes correspondantes.}  $\int_\Delta\, \omega$, ou plus g\'en\'eralement les {\it p\'eriodes-exponentielles}   $\int_\Delta\,  e^{f}\omega$ (o\`u \`a la fois l'int\'egrant $\omega$ et le domaine $\Delta$ sont d\'efinis par des expressions alg\'ebriques d'une ou plusieurs variables \`a coefficients des nombres alg\'ebriques: par exemple  $\pi = \int_0^\infty  \frac{2dt}{1+t^2}$.

  \subsection{}  D'o\`u cette id\'ee vient-elle?

Elle vient d'une th\'eorie initi\'ee par A. Grothendieck (le math\'ematicien qui a r\'evolutionn\'e la g\'eom\'etrie alg\'ebrique dans les ann\'ees 60), la th\'eorie des {\it motifs}, qui vise \`a unifier les aspects combinatoires, topologiques et arithm\'etiques de la g\'eom\'etrie alg\'ebrique. Ces motifs jouent un peu le r\^ole de ``particules \'el\'ementaires" alg\'ebro-g\'eom\'etriques, susceptibles de d\'ecomposition et recombinaison suivant des r\`egles relevant de la th\'eorie des repr\'esentations des groupes. Les groupes en question, baptis\'es {\it groupes de Galois motiviques}, repr\'esentent une formidable {\it g\'en\'eralisation des groupes de Galois usuels aux syst\`emes de plusieurs polyn\^omes \`a plusieurs variables}. Ce ne sont plus des groupes finis, mais des groupes alg\'ebriques (comme dans le cas des groupes de Galois diff\'erentiels).

En tant que groupes de sym\'etrie de motifs (attach\'es \`a de tels syst\`emes \`a coefficients rationnels), ils devraient aussi agir sur leurs p\'eriodes  (et m\^eme sur les avatars exponentiels), ce qui permettrait de d\'efinir les conjugu\'es d'une p\'eriode comme ses images sous les \'el\'ements du groupe de Galois motivique. Toutefois, comme il peut arriver qu'un nombre complexe s'exprime de plusieurs mani\`eres diff\'erentes comme p\'eriode, la coh\'erence d'une telle action requiert de postuler que toute relation entre p\'eriodes provient d'une relation entre motifs, ce qui est la {\it conjecture des p\'eriodes de Grothendieck}. On s'attend en particulier \`a ce que le nombre maximal de p\'eriodes d'un motif sur ${\mathbb Q}$ qui sont alg\'ebriquement ind\'ependantes sur ${\mathbb Q}$ soit \'egal \`a la dimension du groupe de Galois motivique associ\'e\footnote{sur tout cela, on peut consulter l'Introduction aux Motifs de l'auteur, Panoramas et Synth\`eses 17, SMF 2004.}.  

  \subsection{}  Par exemple,  $2\pi\sqrt -1= \int dt/t$ est la p\'eriode attach\'ee au motif de la droite priv\'ee de l'origine, dont le groupe de Galois motivique est le groupe multiplicatif ${\mathbb Q}^\times$, les conjugu\'es de $2\pi\sqrt -1$ \'etant ses multiples rationnels non nuls. Ici, la coh\'erence postul\'ee par la conjecture de Grothendieck n'est autre que la transcendance de $\pi$. 

\medskip Citons aussi \`a titre d'exemples de p\'eriodes, tr\`es \'etudi\'es actuellement mais qui remontent en fait \`a Euler, les nombres polyzeta $$ \sum_{n_1>\dots
>n_k\ge 1}  \,{n_1^{- s_1}\dots n_k^{- s_k}} \; =$$ $$ \; \int_{1\geq t_1 \geq \cdots \geq t_{s} \geq 0}\,  \frac{dt_1}{ \epsilon_1 - t_1}\cdots \frac{dt_s}{ \epsilon_{ s} -  t_{s}}\;\;\;\; (s_i\in  \mathbb Z_{\geq 1},\ ;s=\sum s_i,\; \epsilon_j= 0\; {\text{ou}}\; 1).$$ Depuis Euler, on a d\'ecouvert tout un \'echeveau de relations alg\'ebriques les liant les uns aux autres, et on a v\'erifi\'e que toutes ces relations sont bien d'origine motivique. La th\'eorie motivique sous-jacente aux polyzeta est maintenant bien comprise (gr\^ace surtout aux travaux d'A. Goncharov puis de F. Brown). C'est d'ailleurs le point de vue motivique qui fournit la meilleure majoration $d_s$ - la meilleure connue et conjecturalement la meilleure possible - pour la dimension du ${\mathbb Q}$-espace vectoriel  engendr\'e par  les polyzeta avec $s$ fix\'e: une r\'ecurrence \`a la Fibonacci $d_s= d_{s-2}+d_{s-3}$.

  \subsection{}  Fort heureusement, M. Kontsevich a trouv\'e une formulation \'el\'ementaire, particuli\`erement frappante, de la conjecture des p\'eriodes de Grothendieck, qui ne fait pas appel aux motifs. La voici\footnote{voir aussi M. Kontsevich, D. Zagier: Periods.  Mathematics unlimited---2001 and beyond,  771--808, Springer.}. 

Les deux r\`egles fondamentales du calcul int\'egral que sont la formule de Stokes et le changement (alg\'ebrique) de variables  
 $$\int_\Delta d\omega = \int_{\partial\Delta} \omega,\;\; \;\;\;\int_\Delta f^\ast \omega = \int_{f_\ast \Delta} \omega ,$$  fournissent imm\'ediatement des relations (polyn\^omiales \`a coefficients rationnels) entre p\'eriodes (ou p\'eriodes-exponentielles). La conjecture pr\'edit que, r\'eciproquement, toute relation (polyn\^omiale \`a coefficients rationnels) entre p\'eriodes proviendrait de ces
  deux r\`egles.
  
 \smallskip Par exemple, l'identit\'e d'Euler  $\zeta(2)=\pi^2/6$ peut \^etre comprise comme l'identit\'e de p\'eriodes  $\; \int_0^1\int_0^1\frac{2dxdy}{(1-xy)\sqrt{xy}}= (\int_0^\infty  \frac{2dt}{1+t^2})^2\;$ (par d\'eveloppement en s\'erie g\'eom\'etrique de $1/1-xy$ et int\'egration terme \`a terme, le premier membre s'identifie \`a $6\zeta(2)$). E. Calabi a trouv\'e comment \'etablir cette identit\'e par changement de variables alg\'ebrique idoine: en posant $x=u^2\frac{1+v^2}{1+u^2},\, y= v^2\frac{1+u^2}{1+v^2}$,  de jacobien $\vert \frac{d(x,y)}{d(u,v)}\vert = \frac{4uv(1-u^2v^2)}{(1+u^2)(1+v^2)}= \frac{4(1-xy)\sqrt{xy}}{(1+u^2)(1+v^2)}$, on obtient  $ \int_0^1\int_0^1\frac{2dxdy}{(1-xy)\sqrt{xy}}=   \int\int_{u,v\geq 0, uv\leq 1}\, \frac{8dudv}{(1+u^2)(1+v^2)}$,  qu'on identifie \`a  $(\int_0^\infty  \frac{2dt}{1+t^2})^2$  en utilisant l'involution $(u,v\mapsto u^{-1}, v^{-1}$)\footnote{C. Viola vient de me communiquer une variante encore plus simple, consistant \`a identifier $\zeta(2)$ \`a $ \int_0^1\int_0^1\frac{dxdy}{(1-xy}$, et poser $x=u-v, y=u+v$, ce qui donne $\zeta(2)/4= \int_0^{1/2} du \int_0^u \frac{dv}{1-u^2+v^2}+   \int_{1/2}^{1} du \int_0^{1-u} \frac{dv}{1-u^2+v^2}$, qu'on \'evalue par des moyens \'el\'ementaires \`a $\pi^2/24$.}.

   \subsection{}  En r\'esum\'e, on arrive \`a cette id\'ee sp\'eculative g\'en\'erale, tout \`a fait dans l'esprit de la lettre-testament de Galois, que {\it l'arithm\'etique de cette vaste classe de nombres - les p\'eriodes-exponentielles - devrait \^etre dict\'ee par les r\`egles \'el\'ementaires du calcul int\'egral (et s'il en est ainsi, se d\'ecrire en termes de groupes {\og}{galoisiens}{\fg})}.

 \smallskip On peut aller un peu plus loin, comme l'a montr\'e tout r\'ecemment J. Ayoub, en ne retenant que la r\`egle de Stokes: la r\`egle du changement de variable en d\'ecoule. 

Avec d'expliquer pourquoi, voici la formulation qu'Ayoub donne de la conjecture des p\'eriodes.
Soit $\mathcal A$ l'alg\`ebre des fonctions d'un nombre quelconque de variables complexes $z_i$, qui sont holomorphes dans le polydisque $z_i\leq 1$, et alg\'ebriques sur $\mathbb Q(z_1, \ldots, z_i,\ldots)$. 
On a une forme lin\'eaire 
$\; \mathcal A  \stackrel{\int_{\square}}\longrightarrow   {\mathbb C} \;$ 
donn\'ee par l'int\'egrale sur l'hypercube r\'eel $z_i\in [0,1]$, et il s'agit de d\'ecrire son noyau.

Pour tout indice $i$ et toute fonction $g_i\in \mathcal A$, la fonction $h_i = \frac{\partial g_i}{\partial z_i} - g_{i\mid z_i= 1}+ g_{i\mid z_i= 0}$ est dans le noyau de $\int_{\square}$. La conjecture pr\'edit que, r\'eciproquement, tout \'el\'ement du noyau de $\int_{\square}$ serait combinaison lin\'eaire de telles $h_i$. 

 \smallskip Par exemple, partant de $h(z_1)\in \mathcal A$ et d'une fonction $f(z_1)$ alg\'ebrique qui envoie le dique unit\'e (\resp l'intervalle $[0,1]$) dans lui-m\^eme en fixant $0$ et $1$, la formule de changement de variable montre que  $ f'(z_1)h(f(z_1)) - h(z_1)$ est dans le noyau de $\int_{\square}$. 
 
 Comme l'a observ\'e Ayoub, on peut, quitte \`a ajouter une variable, \'ecrire $ f'(z_1)h(z_1) - h'(z_1)$ sous la forme pr\'edite par la conjecture: poser $f_1= f(z_1)-z_1,\; f_2= -z_2f'(z_1)+ z_2-1,$ et $g_i = f_i\cdot h(z_2f(z_1)+(1-z_2)z_1)$ pour $i=1,2$. 
  
 \section{Coda: un groupe de Galois {\og}cosmique{\fg}?} Depuis quelques temps, les id\'ees galoisiennes ont fait irruption en physique quantique, plus pr\'ecis\'ement en th\'eorie perturbative des champs quantiques. 
 
 La divergence, plus importune encore en physique qu'en math\'ematique, infeste la physique des champs quantiques. \`A partir des travaux de R. Feynman et J. Schwinger,  les physiciens ont d\^u b\^atir un arsenal de techniques, bien plus \'elabor\'ees que celles de Euler, pour la d\'epasser.  Le proc\'ed\'e le plus simple et le plus utilis\'e est la renormalisation par r\'egularisation dimensionnelle:  on fait fluctuer la dimension de l'espace-temps en lui faisant prendre des valeurs complexes voisines de $4$, et on d\'eveloppe les int\'egrales obtenues en s\'eries index\'ees par des diagrammes de Feynman de complexit\'e croissante. L'\'elimination des termes {\og}{divergents}{\fg} de la s\'erie se fait suivant de subtiles r\`egles combinatoires qui garantissent la coh\'erence du proc\'ed\'e. 

  En jonglant avec des int\'egrales divergentes, donc d\'epourvues de sens physique (et m\^eme math\'ematique, {\it a priori}), la renormalisation aboutit \`a des quantit\'es finies, en accord remarquable avec l'exp\'erience de surcro\^{\i}t. Elle r\'eussit le tour de force d'{\it extraire, syst\'ematiquement, du (d\'e)fini de l'in(d\'e)fini}.

\smallskip La {\og}{m\oe lle}{\fg} math\'ematique de ce proc\'ed\'e a r\'ecemment \'et\'e extraite par A. Connes et D. Kreimer, qui ont associ\'e aux th\'eories quantiques des champs des groupes de sym\'etries infinis directement construits en termes de diagrammes de Feynman. En effectuant une {\og}{mont\'ee vers l'absolu}{\fg}, ils obtiennent un groupe de Galois absolu - le groupe de Galois {\og}{cosmique}{\fg} pr\'edit par P. Cartier\footnote{{\og}{La folle journ\'ee, de Grothendieck \`a Connes et Kontsevich}{\fg}, Festschrift des 40 ans de l'IHES.} - qui agit sur les {\og}{constantes}{\fg} de toutes les th\'eories quantiques des champs. 

Ce groupe, d'une ubiquit\'e stup\'efiante,  incarne \`a lui seul les divers avatars galoisiens \'evoqu\'es ci-dessus:

- il s'interpr\`ete comme groupe de Galois diff\'erentiel.

- il est tr\`es proche du groupe de Galois motivique attach\'e aux nombres polyz\^etas (nombres qu'on retrouve souvent en calculant des int\'egrales de Feynman).

 - c'est une variante alg\'ebro-g\'eom\'etrique du groupe ${\rm{GT}}$ de Drinfeld.

 \medskip C'est ainsi qu'en th\'eorie quantique des champs, les divergences, loin d'\^etre des nuisances, donnent naissance \`a des ambigu\"{\i}t\'es galoisiennes formant le groupe de sym\'etries d'une riche structure qui appara\^{\i}t dans des domaines math\'ematiques tr\`es \'eloign\'es les uns des autres.

     \end{document}